\documentclass[a4paper,12pt]{article}

\usepackage{listings}
\usepackage{longtable}

\usepackage{graphicx}

\usepackage[font=scriptsize]{caption}

\usepackage{float}

\usepackage{amsmath,amssymb,mathrsfs,eucal}

\usepackage[T1,T2A]{fontenc}

\font\cyrfont=wncyss10

\def\sza{\hbox{\cyrfont X}} 

\newtheorem{thm}{Theorem}

\newtheorem{conj}[thm]{Conjecture}

\usepackage{color}

\begin{document}

\title{Critical $L$-values for some quadratic twists of Gross curves}

\author{Andrzej D\k{a}browski, Tomasz J\k{e}drzejak and Lucjan Szymaszkiewicz} 

\date{}

\maketitle{}

Let $K=\Bbb Q(\sqrt{-q})$, where $q$ is a prime congruent to $3$ modulo $4$. Let $A=A(q)$ 
denote the Gross curve \cite{Gr}. Let $E=A^{(-\beta)}$ denote its quadratic twist, with 
$\beta=\sqrt{-q}$. The curve $E$ has the nice explicit equation (see  \cite{CL}, equation (1.2)) 

\begin{equation}\label{e:01} 
y^2 = x^3 - 2^{-4}3^{-1}(j({\mathcal O}_K)^{1/3}x + 2^{-5}3^{-3}(j({\mathcal O}_K) - (12)^3)^{1/2}, 
\end{equation} 
where it is understood that, in this equation, we take the real cube root of $j({\mathcal O}_K)$, 
and the square root of $j({\mathcal O}_K) - (12)^3$ lying in the upper half complex plane. Thus $E$ 
is defined over the Hilbert class field $H$ of $K$. Below we use Magma \cite{BCP} 
to calculate the values $L(E/H,1)$ for all such $q$'s up to some reasonable ranges 
(different for $q\equiv 7 \, \text{mod} \, 8$ and  $q\equiv 3 \, \text{mod} \, 8$). 
All these values are non-zero, and using the Birch and Swinnerton-Dyer conjecture, 
we can calculate hypothetical orders of $\sza(E/H)$ in these cases. Our calculations extend those 
given in \cite{CC} for the case $q=7$.

\section{The case $q\equiv 7 \, \text{mod} \, 8$}

In this case we know, by a recent result of J. Coates and Y. Li (\cite{CL}, Theorem 1.3), 
that $L(E/H,1) \not=0$. In the table below we calculate numerically these values for all such $q$ 
up to $4663$.  

Now let us say a few words about the Magma implementation. The starting source for us 
was the article by M. Watkins \cite{Wa}, which gives some numerical examples 
(or rather hints) how to compute 
Grossencharacters and critical L-values (sections 5.4 and 6.1 deals with 
$\mathbb Q(\sqrt{-23})$, but of course we need to keep track of the effect of twisting). 
But it was not enough for us to write an algorithm calculating $L(E/H,1)$. Watkins \cite{Wat}
corrected our algorithm (or better, he wrote a new one) and 
tested for $q=23$ and $79$. It was a starting point for us to make extensive 
numerical calculations. The algorithm uses the fact that 
 $L$-series of an elliptic curve over $H$ splits into factors 
corresponding to Grossencharacters twisted by Hilbert 
characters and its conjugates (so it uses the classical Hecke-Deuring theory 
linking elliptic curves with $CM$ to Grossencharacters, with keeping 
track of the effect of twisting). Here are some more details. 
Assuming $2{\mathcal O}_K = {\mathfrak p}{\mathfrak p}^{\star}$, 
and choosing the sign of $\sqrt{-q}$ so that $\text{ord}_{\mathfrak p}((1-\beta))>0$, 
we can check that $E/H$ has good reduction outside the primes of $H$ 
lying above $\mathfrak p$ (see \cite{CL}). Moreover, the Deuring Grossencharacter 
$\psi_{E/H}$ of $E/H$ is then equal to $\rho \circ N_{H/K}$, where 
$\rho$ is the Grossencharacter of $K$ with conductor $\mathfrak p^2$ defined by 
$$
\rho(\mathfrak a) = \alpha, \quad \mathfrak a^h = \alpha\mathcal O_K, \quad 
\alpha \equiv 1 \, \text{mod} \, \mathfrak p^2. 
$$
The algorithm computes the values $L(\rho\chi,1)$, where $\chi$ runs over the characters 
of the ideal class group of $K$. Now thanks to the above formula and Deuring's 
theory it follows that $L(E/H,1)$ will be given by the product of all the $L(\rho\chi,1)$'s  
and their complex conjugates.

\bigskip

Now we know by Iwasawa theory that the Tate-Shafarevich group $\sza(E/H)$ is finite because 
$L(E/H,1) \not=0$. Below we will write down an explicit conjectural formula for the order of 
$\sza(E/H)$.  Let $h$ denote the class number of $K$, $m=\frac{q-1}4 - \frac h 2$, and let  
$\Omega(q)$ be a period defined in \cite{Gr}: 

$$
\Omega(q) = \frac 1 { (2 \pi)^m \cdot q^{\frac h 2} }
\prod_{ \substack{0<c<q \\ \left(\frac c q\right)=1} }
\Gamma \left(\frac c q \right). 
$$ 
The prime $2$ splits in $K$, and we write $2{\mathcal O}_K = {\mathfrak p}{\mathfrak p}^{\star}$, 
where we have chosen the sign of $\beta=\sqrt{-q}$ so that $\text{ord}_{\mathfrak p}((1-\beta))>0$. 
Then $E$ is the quadratic twist of the Gross curve $A/H$ by $H(\sqrt{-\beta})/H$. 
Let $\{v_1,...,v_r\}$ be the set of primes of $H$ lying above $\mathfrak p$, so that 
$r=h/j$, where $j$ is the exact order of the class of $\mathfrak p$ in the ideal class group 
of $K$.  It turned out that there are exactly $18$ primes $q \leq 4663$ congruent to $7$ 
modulo $8$, for which $r > 1$.

\begin{conj} 
$\#(\sza(E/H)) = L(E/H,1) 2^{h+6-2r} / (\Omega(q)^2\sqrt{q})$. 
\end{conj} 

One easily checks that in case $q=7$, the formula from Conjecture 1 is equivalent to 
(2.11) from \cite{CC}. 
The above conjecture agrees with the Birch and Swinnerton-Dyer conjecture for $E/H$. 
In particular, the set of bad primes of $E/H$ is precisely $\{v_1,...,v_r\}$, and the 
factor $2^{-2r}$ in the above formula takes account of the fact that the Tamagawa 
factor at each of these primes is $4$. 

\bigskip 
 
John Coates informed one of us (A. D.) that one should be able to use the Iwasawa 
theory being developed in \cite{CKLT} to prove the above conjecture.

\bigskip  

Below we use Conjecture 1 to calculate $\#(\sza(E/H))$ (i.e. the analytic order of 
$\sza(E/H)$) for all primes $q$ congruent to $7$ modulo $8$ up to $4663$.

\begin{center}
\tiny
\begin{longtable}{|r|r|r|r|} 
\hline 
\multicolumn{1}{|c|}{$q$} 
& 
\multicolumn{1}{|c|}{$h$} 
&
\multicolumn{1}{|c|}{$L(E/H,1)$} 
& 
\multicolumn{1}{|c|}{$\sqrt{\#(\sza(E/H))}$} 
\\ 
\hline 
\endhead 
\hline 
\multicolumn{4}{|r|}{{Continued on next page}} \\
\hline
\endfoot

\hline 
\hline
\endlastfoot

7 &	1 &	0.30903153751765917103 &	1 \\
23 &	3 &	0.79196294535428296044 &	1 \\
31 &	3 &	0.35288571505654851763 &	1 \\
47 &	5 &	3.25049251883301426121 &	3 \\
71 &	7 &	0.10920125590289049507 &	1 \\
79 &	5 &	0.02577591231345318312 &	1 \\
103 &	5 &	0.84244014254446144514 &	13 \\
127 &	5 &	0.33138747507581642444 &	17 \\
151 &	7 &	0.00899919291175804982 &	5 \\
167 &	11 &	338.84342541058916626822 &	2049 \\
191 &	13 &	0.07538843930533773444 &	81 \\
199 &	9 &	0.00178784908116291475 &	9 \\
223 &	7 &	0.66858391145992740299 &	289 \\
239 &	15 &	0.02401256252449269664 &	311 \\
263 &	13 &	0.24799355777337639904 &	1767 \\
271 &	11 &	0.00495300516895988511 &	127 \\
311 &	19 &	0.08289319536914465106 &	12559 \\
359 &	19 &	0.00008262935013341212 &	2057 \\
367 &	9 &	0.02393861560648477609 &	1679 \\
383 &	17 &	1058.78512825720370837609 &	9090067 \\
431 &	21 &	5.38876192180481196261 &	2039928 \\
439 &	15 &	0.00002907103487183395 &	1279 \\
463 &	7 &	0.08500423491817571054 &	11663 \\
479 &	25 &	1483.07868786791841546796 &	1746287691 \\ 
487 &	7 &	0.14694723669623207042 &	22807 \\
503 &	21 &	260759.24737728583571680044 &	2880463783 \\ 
599 &	25 &	0.00000001076991986883 &	162285 \\
607 &	13 &	0.10795424186869623536 &	884605 \\
631 &	13 &	0.00004385443164140780 &	44425 \\
647 &	23 &	0.00000607641351529086 &	4925391 \\
719 &	31 &	0.00530561250904147645 &	31646320057 \\
727 &	13 &	0.00299892234779012135 &	1113693 \\
743 &	21 &	1498.05565627050935641062 &	332146468299 \\
751 &	15 &	0.00000000896688802629 &	10512 \\
823 &	9 &	0.00102716234866514602 &	469855 \\
839 &	33 &	0.12347121525795507868 &	9315485111867 \\
863 &	21 &	4.28771850132666368851 &	240267975371 \\
887 &	29 &	19477038.03896518654448808873 &	30785347392739103 \\
911 &	31 &	0.00000312260876667640 &	98895319091 \\
919 &	19 &	0.00000000001480555834 &	36741 \\
967 &	11 &	0.09804744238584120050 &	67762715 \\
983 &	27 &	80292384.19994893868777292158 &	238138622744502833 \\
991 &	17 &	0.00004920969572684842 &	91561037 \\
1031 &	35 &	0.06252031110591168794 &	1595084268489133 \\
1039 &	23 &	0.00000369632068956271 &	1196802971 \\
1063 &	19 &	1.67498940035011217030 &	54030361471 \\
1087 &	9 &	0.00390824946108466869 &	22866381 \\
1103 &	23 &	0.00000572944338432493 &	80104082513 \\
1151 &	41 &	0.16613607401358903141 &	1320694330164429335 \\
1223 &	35 &	0.00194039562173745137 &	9307326019643999 \\
1231 &	27 &	0.00000000093497039710 &	3758406353 \\
1279 &	23 &	0.00000000287686176036 &	1814619877 \\
1303 &	11 &	0.00021448455521015721 &	157522307 \\
1319 &	45 &	28941477.30398577348268950039 &	494146202273056454638285 \\ 
1327 &	15 &	0.00028916732527519815 &	1982319913 \\
1367 &	25 &	28601.01936469536020713962 &	1839953001047559952 \\
1399 &	27 &	0.00000001208663958591 &	210925423356 \\
1423 &	9 &	0.01738234231379547068 &	1715895524 \\
1439 &	39 &	0.00000006661042814296 &	30644209290657623 \\ 
1447 &	23 &	0.00013884708802240947 &	811815992737 \\
1471 &	23 &	0.00000000029882112555 &	5374292551 \\
1487 &	37 &	6576514753.03121783833652353519 &	10937565775616401748256581 \\
1511 &	49 &	0.00055060129228412947 &	25127641761490803096975 \\
1543 &	19 &	0.00176400518665963757 &	1034981247929 \\
1559 &	51 &	0.87726903612746376846 &	2738379667079823194690949 \\ 
1567 &	15 &	0.20108355789279651145 &	1575131870837 \\
1583 &	33 &	0.13269489544044152769 &	15757834130321482863 \\
1607 &	27 &	0.55194730530618115640 &	349749435110817399 \\ 
1663 &	17 &	0.00119021469463015769 &	513229201123 \\
1759 &	27 &	0.00000000101913474126 &	1344751145077 \\ 
1783 &	17 &	0.00024486558886747311 &	1885276250643 \\
1823 &	45 &	9746102.85279061269213579273 &	1808382408390942813565214784 \\ 
1831 &	19 &	0.00000000000009270267 &	880904745 \\
1847 &	43 &	1207965.54129407491879223532 &	4921203932360352431818054853 \\
1871 &	45 &	0.00000467957783395582 &	211817576664143143791049 \\
1879 &	27 &	0.00000055405311625279 &	359706350854836 \\ 
1951 &	33 &	0.00000000000191765783 &	71833637298811 \\ 
1999 &	27 &	0.00000000047804213103 &	21632373948999 \\ 
2039 &	45 &	0.00000000367413120950 &	7999274603520740597625 \\ 
2063 &	45 &	2766645.80955923764536501056 &	1341376803998421383216338486800 \\
2087 &	35 &	164.98809160497598617071 &	7549276579545660794505221 \\
2111 &	49 &	0.00109168602722010353 &	4327178311989716140475456427 \\
2143 &	13 &	0.01476049903300890724 &	12322652032019 \\
2207 &	39 &	8356.57170711443065175664 &	5015914871462321771139549955 \\
2239 &	35 &	0.00000000000004393777 &	11915930073079139 \\
2287 &	29 &	0.00766729713419227393 &	6824315445622624657 \\
2311 &	29 &	0.00000000187662881760 &	20015174351579955 \\
2351 &	63 &	552470425.54422669425872718713 &	2816977737852577527309898650097367443659 \\ 
2383 &	29 &	0.00002630161314175636 &	1162866812268136121 \\
2399 &	59 &	0.04141467915645522162 &	13235014758527507066087807686368973 \\
2423 &	33 &	250.16970151756207630008 &	281045797696507991427800573 \\
2447 &	37 &	23.82978667301734246096 &	2092852385170614860214887843 \\
2503 &	21 &	0.02648336374331768867 &	100742172570243491 \\
2543 &	35 &	382.95959917323265302185 &	10157978602722623458338530496 \\
2551 &	41 &	0.00000000000008248960 &	60637243039930628065 \\
2591 &	57 &	0.00000000000238017737 &	336275233499026658311829296679 \\
2647 &	15 &	0.03372248437933683851 &	6798505360748663 \\
2663 &	43 &	822.32375569716777860250 &	76481720773755781866561562664447 \\
2671 &	23 &	0.00000000000003162935 &	66697260897735 \\
2687 &	51 &	194688424.12336279094797547720 &	17321740551070983627343715625416872652 \\
2711 &	53 &	0.37421143332764690649 &	38629854923911845038305204026873649 \\
2719 &	41 &	0.00000000000083400226 &	1406443138529321307393 \\
2767 &	21 &	0.31592061458094422698 &	3800789493866165232 \\
2791 &	39 &	0.00000000000054815909 &	897017319101036106209 \\
2879 &	57 &	0.00000013893708474299 &	10273597469646245008601935022410689 \\
2887 &	25 &	0.00000727449204224937 &	4682281660959493201 \\
2903 &	59 &	469211113283732.46988840553969223592 &	400977393247118763374214959073343661518898573 \\
2927 &	31 &	64.29980720973840221765 &	27260274502712067357277082151 \\
2999 &	73 &	360.95151202061450399689 &	3943185992249268695560714545947130529724467545 \\
3023 &	47 &	511953545751.66177107863339231208 &	6789457020754215256411685440174304374421 \\
3079 &	41 &	0.00000000023207164115 &	2141735462113012348434191 \\
3119 &	69 &	0.00602643768003535676 &	1049324059291363104640660260593606354947228 \\
3167 &	53 &	260425057231369.09972689242486572823 &	179940200099564565183145128420460525695159379 \\
3191 &	69 &	14.95753417523273502334 &	628315606611285122323089495662377335179881369 \\
3271 &	27 &	0.00000001241846857681 &	243743805808037264960 \\
3319 &	41 &	0.00000000000000931719 &	33642982599580090033407 \\
3343 &	19 &	1.58620548101554534230 &	202212319807498314139 \\
3359 &	69 &	1.67365772978370467795 &	4049500549734207136741889632539175584626802773 \\
3391 &	37 &	0.00000000047269143915 &	3101068731835541701552139 \\
3407 &	57 &	28819363.52788452185843622630 &	14560069825940266325266106631738472810251876 \\ 
3463 &	19 &	13.77281385527310381528 &	8436425008239920177883 \\
3511 &	41 &	0.00000000001434615685 &	11283089506997160477570507 \\
3527 &	65 &	7976156006.93765520950738207851 &	177909675370561108532941247559470051694366241471 \\ 
3559 &	45 &	0.00000000000454262476 &	2360984303797804985269627200 \\
3583 &	29 &	0.03920861329093863272 &	5633558142934297942544477 \\
3607 &	19 &	36.90684229826591202567 &	6563789114878102952373 \\
3623 &	45 &	26.02198199101404503659 &	1367770373787759276795854320370634923 \\ 
3631 &	43 &	0.00000000000255623850 &	427172018608825201653180173 \\
3671 &	81 &	68467807.61967835918310312766 &	117843219225119087861773569889261690958967498129448976645 \\
3719 &	67 &	0.00000000387780125945 &	2100976923197097618349171255793507783224063 \\
3727 &	31 &	0.00004718946253999069 &	2444912707864422958673293 \\
3767 &	39 &	0.14327587591455739687 &	19175435293484578149783919131747899 \\
3823 &	29 &	0.00002507286357856734 &	112417612198200717861969 \\
3847 &	23 &	0.00001533345082846095 &	862130992415445336857 \\
3863 &	61 &	15711737.48318757316819605069 &	538394367773173591693084396367727593976141647921 \\
3911 &	83 &	0.00000000001135189188 &	14499327197980399167915624859957367185588546363629 \\
3919 &	39 &	0.00000000044261630153 &	3669196191866220345463135468 \\
3943 &	27 &	0.00000000002095427812 &	24733668567484868147 \\
3967 &	33 &	159215.22350765659455009473 &	1699553443719169549240379327572 \\ 
4007 &	57 &	0.91990088235951608333 &	24479483802292931184510605398988563448578689 \\
4079 &	85 &	0.00000127446973332763 &	321963672249515003702405195900159942517697859377925616 \\
4111 &	39 &	0.00000000000002095214 &	33609809857361575194185877 \\
4127 &	49 &	264235.35416434781417718265 &	8507300729794243618442185373646968515974841 \\
4159 &	31 &	0.00000000108677317318 &	8139478924692126488523677 \\
4231 &	51 &	0.00000000799222631822 &	4299570633402618922705497439694423 \\ 
4271 &	65 &	0.00000290128128404136 &	28283186967833318171568219977944584501489867447 \\
4327 &	19 &	0.00017450275228592056 &	1066852592208273415311 \\
4391 &	79 &	0.00001251482541320985 &	111218267845699128837861765031864433022917071402312905 \\
4423 &	33 &	1.13324073727367674386 &	572395522360267105755996447652 \\
4447 &	17 &	0.00000077492121857457 &	39872685366747226231 \\
4463 &	55 &	115564.52263908769265426894 &	72780887497608197802303423424278364999654385253 \\
4519 &	29 &	0.00000000000000004313 &	3597009679993991314033 \\
4567 &	33 &	0.00000017729987000499 &	1945455740103468168767260495 \\
4583 &	61 &	551441731499.78356036276104299918 &	1438973413788257170719455133810181008462888610212680883 \\
4591 &	49 &	0.00000000000000079237 &	16108217968515978652127771698061 \\
4639 &	51 &	0.00000000000000004301 &	101877348949955205680678906825472 \\
4663 &	33 &	0.00001145384710376496 &	1169504442816257396334162500 \\

\end{longtable}
\end{center}

\section{The case $q \equiv 3 \, \text{mod} \, 8$}

In this case the curve $E$ defined by the equation (1) is also defined over $H$. Here the prime  
$2$ is inert in $K$, and the curve $E$ will always have good reduction outside the set of 
the primes of $H$ lying above $2$  (assuming $q>3$).  The Deuring Grossencharacter 
$\psi_{E/H}$ of $E/H$ is then equal to $\rho \circ N_{H/K}$, where 
$\rho$ is the Grossencharacter of $K$ with conductor $4\mathcal O_K$ defined by 
$$
\rho(\mathfrak a) = \alpha, \quad \mathfrak a^h = \alpha\mathcal O_K, \quad 
\alpha \equiv 1 \, \text{mod} \, 4\mathcal O_K. 
$$
The algorithm computes the values $L(\rho\chi,1)$, where $\chi$ runs over the characters 
of the ideal class group of $K$. Again, thanks to the above formula and Deuring's 
theory it follows that $L(E/H,1)$ will be given by the product of all the $L(\rho\chi,1)$'s  
and their complex conjugates. 

Our numerical calculations (given in the table below) lead to the following conjecture 
(see \cite{CL}, Conjecture 1.5). 

\begin{conj} 
For all primes $q$ with $q\equiv 3 \, \text{mod} \, 8$, we have $L(E/H,1)\not=0$. 
\end{conj}

As it is remarked in (\cite{CL}, p. 2), in contrast to the proof of Theorem 1.3 there, the 
authors see no way at present for attacking such a conjecture using Iwasawa theory.  

In this case, we propose the following conjectural formula for the order of $\sza(E/H)$. 
Now, the Tamagawa factor at each prime of bad reduction is $1$, and the following 
conjecture agrees with the Birch and Swinnerton-Dyer conjecture for $E/H$.  

\begin{conj} 
$\#(\sza(E/H)) = L(E/H,1) 2^{2h} / (\Omega(q)^2\sqrt{q})$. 
\end{conj}

Below we use Conjecture 3 to calculate $\#(\sza(E/H))$ (i.e. the analytic order of 
$\sza(E/H)$) for all primes $q$ congruent to $3$ modulo $8$ up to $11131$.

\begin{center}
\tiny
\begin{longtable}{|r|r|r|r|} 
\hline 
\multicolumn{1}{|c|}{$q$} 
& 
\multicolumn{1}{|c|}{$h$} 
&
\multicolumn{1}{|c|}{$L(E/H,1)$} 
& 
\multicolumn{1}{|c|}{$\sqrt{\#(\sza(E/H))}$} 
\\ 
\hline 
\endhead 
\hline 
\multicolumn{4}{|r|}{{Continued on next page}} \\
\hline
\endfoot

\hline 
\hline
\endlastfoot 

11 &	1 &	1.73845792121760807790 &	1 \\
19 &	1 &	1.00717576250064706853 &	1 \\
43 &	1 &	1.27416027648354776885 &	2 \\
59 &	3 &	3.27291981598555587930 &	5 \\
67 &	1 &	1.22243364144817892444 &	3 \\
83 &	3 &	0.03764760689032642372 &	1 \\
107 &	3 &	1.03936693115122887083 &	9 \\
131 &	5 &	0.00697158270425921776 &	2 \\
139 &	3 &	0.06083483360363034315 &	3 \\
163 &	1 &	1.23339224060989144293 &	10 \\
179 &	5 &	0.06195226194655516123 &	17 \\
211 &	3 &	0.00822152771334987023 &	3 \\
227 &	5 &	101.51957725718817748183 &	1524 \\
251 &	7 &	0.32436912039697336340 &	343 \\
283 &	3 &	0.66685505848109321080 &	53 \\
307 &	3 &	0.78609012909015461469 &	73 \\
331 &	3 &	0.14907707956876846359 &	49 \\
347 &	5 &	0.08080910761305406985 &	250 \\
379 &	3 &	0.33390038219003027566 &	117 \\
419 &	9 &	1.67741166405082133222 &	40225 \\
443 &	5 &	0.00285060561331757659 &	156 \\
467 &	7 &	0.35140804797178265726 &	8190 \\
491 &	9 &	0.01884246766690787210 &	11387 \\
499 &	3 &	0.46414182882437918271 &	395 \\
523 &	5 &	0.02170338307342622816 &	242 \\
547 &	3 &	0.41073019842431356261 &	316 \\
563 &	9 &	2217.38322210227025102910 &	7489893 \\
571 &	5 &	0.29829944260019244950 &	1956 \\
587 &	7 &	0.25898213029907196441 &	31143 \\
619 &	5 &	0.00509394359043077170 &	479 \\
643 &	3 &	0.60354061296945788027 &	893 \\
659 &	11 &	1.79783553470780040531 &	4006249 \\
683 &	5 &	0.00159958302432158035 &	1447 \\
691 &	5 &	0.10350577224179846609 &	4018 \\
739 &	5 &	0.10706905313737192318 &	4561 \\
787 &	5 &	0.01395924800465447269 &	1425 \\
811 &	7 &	0.17939743855163310334 &	60156 \\
827 &	7 &	154.81828584920445528361 &	7220288 \\
859 &	7 &	0.00003489503731419296 &	1202 \\
883 &	3 &	1.46744585216089836663 &	4031 \\
907 &	3 &	0.59650901915093688786 &	2888 \\
947 &	5 &	0.03074953029629406725 &	62879 \\
971 &	15 &	0.04674641781017625752 &	846868715 \\
1019 &	13 &	0.23524712311785247652 &	481254336 \\
1051 &	5 &	0.15568884564868710538 &	48000 \\
1091 &	17 &	0.00217725207751161362 &	4271088999 \\
1123 &	5 &	0.07364280939031303020 &	23322 \\
1163 &	7 &	0.37473402025732946137 &	6435488 \\
1171 &	7 &	0.00082674907382398778 &	33427 \\
1187 &	9 &	8.57806631683350786688 &	196575884 \\
1259 &	15 &	0.04793034848358401730 &	23807653018 \\
1283 &	11 &	8.90099754034007899869 &	2309889447 \\
1291 &	9 &	0.00241815361075660441 &	1270411 \\
1307 &	11 &	16150.89301909056086676059 &	186585360146 \\
1427 &	15 &	1.56491207854150961696 &	208776957205 \\
1451 &	13 &	5.78796172237250084560 &	198287772553 \\
1459 &	11 &	0.02184103208570492546 &	81647642 \\
1483 &	7 &	3.19082650762232731055 &	6018956 \\
1499 &	13 &	0.00107052328295353035 &	3762443612 \\
1523 &	7 &	13.14472327791761386093 &	432126977 \\
1531 &	11 &	0.00934308869252333838 &	91517708 \\
1571 &	17 &	77374.81778094942048006157 &	4227257131159937 \\
1579 &	9 &	0.04967065171792227175 &	18842731 \\
1619 &	15 &	0.01588367050528775720 &	287671778919 \\
1627 &	7 &	1.14234502669777570952 &	8445329 \\
1667 &	13 &	3015.72287790923136695989 &	5415661616355 \\
1699 &	11 &	0.00000951105549164743 &	7161466 \\
1723 &	5 &	0.59636225379566621680 &	627691 \\
1747 &	5 &	0.63250816243545283913 &	545202 \\
1787 &	7 &	0.85334738271374115107 &	618910949 \\
1811 &	23 &	0.06992899432395535990 &	16061273092160342 \\
1867 &	5 &	2.38591921224867532372 &	6297927 \\
1907 &	13 &	6.99689763033106016118 &	3165747533075 \\
1931 &	21 &	90.38028537368310596632 &	258729934559477369 \\
1979 &	23 &	4142.79431454502235459380 &	21986676642147049263 \\
1987 &	7 &	0.00236988066983797578 &	2280430 \\
2003 &	9 &	0.01742633984133732059 &	2517289989 \\
2011 &	7 &	0.00053306208627154529 &	2638203 \\
2027 &	11 &	0.14670352596641644330 &	58540694152 \\
2083 &	7 &	0.05598543472754948453 &	4585113 \\
2099 &	19 &	0.03659270303763838186 &	2787484250243039 \\
2131 &	13 &	0.08176766229504099310 &	89178286597 \\
2179 &	7 &	0.00000665857111220152 &	527657 \\
2203 &	5 &	1.91715188964021586919 &	15965226 \\
2243 &	15 &	455169.94065435589748159377 &	81753149480991824 \\
2251 &	7 &	0.02300309205779565400 &	64040332 \\
2267 &	11 &	7.28105312585640396026 &	4870756173147 \\
2339 &	19 &	0.97878697785561566356 &	59954610410365278 \\
2347 &	5 &	1.60476528633349239697 &	7286700 \\
2371 &	13 &	0.00000800495831506833 &	3115285126 \\
2411 &	23 &	504.85388397709402086354 &	362000857310467738783 \\
2459 &	19 &	0.93923549011641697155 &	117489182977955129 \\
2467 &	7 &	0.09990402861216668317 &	100745543 \\
2531 &	17 &	0.00014747927341115960 &	306040338160927 \\
2539 &	11 &	0.00011234010687514967 &	3085577520 \\
2579 &	21 &	15539.25782562482093805923 &	758227028237459575815 \\
2659 &	13 &	0.14349748145910982056 &	1922642144666 \\
2683 &	5 &	0.35221319600572141285 &	14543262 \\
2699 &	15 &	0.00003740841184509724 &	50729716969650 \\
2707 &	7 &	0.17461587946412744900 &	348539119 \\
2731 &	11 &	0.00000353700779516742 &	321349458 \\
2803 &	9 &	0.78284831007160905719 &	6224866339 \\
2819 &	21 &	0.27671024228777928761 &	20477277172589698869 \\
2843 &	15 &	10.14431840386199781185 &	4022140478178599 \\
2851 &	11 &	0.00056503632880051287 &	11870467976 \\
2939 &	29 &	0.00624713894355604616 &	105268137875003312953547 \\
2963 &	13 &	1.32077115136644560933 &	319550765817053 \\
2971 &	11 &	0.03190494829015117475 &	324491996326 \\
3011 &	21 &	0.26180404016133656775 &	33950238559165464896 \\
3019 &	7 &	0.04312030585447679517 &	1802300809 \\
3067 &	7 &	0.50159638459199363857 &	383973387 \\
3083 &	13 &	66.44365839794315801657 &	15879087340845386 \\
3163 &	9 &	0.00000102202290136419 &	76089107 \\
3187 &	7 &	0.21740590516706865825 &	569775265 \\
3203 &	11 &	7.58674726638825494993 &	378391105040096 \\
3251 &	31 &	1.49283692874036018378 &	208144942072228395025506250 \\
3259 &	9 &	0.01121435077834574673 &	21243653932 \\
3299 &	27 &	0.00000054505993369606 &	403911133854039617472 \\
3307 &	9 &	5.24935939959546086799 &	190357002279 \\
3323 &	17 &	7353.12213893642943109960 &	56954453746238785334 \\
3331 &	15 &	0.00013628808312750543 &	12824146340774 \\
3347 &	11 &	0.00693687813140206147 &	15997504388590 \\
3371 &	21 &	0.00000905587936280291 &	1798520024157799359 \\
3467 &	19 &	0.14118916689423367305 &	2179702892746704617 \\
3491 &	23 &	1521.81323914654475213527 &	2216917282043791853812677 \\
3499 &	11 &	0.00000048245603947392 &	12562962745 \\
3539 &	23 &	0.00000000024192865064 &	502868601762633637 \\
3547 &	9 &	0.51723990155998270507 &	64134211184 \\
3571 &	15 &	0.00001877157784780692 &	5298557564473 \\
3643 &	9 &	0.19035162141447944342 &	196656106527 \\
3659 &	29 &	0.52058334069619626047 &	224345918498470661616572059 \\
3691 &	13 &	0.00010142020688338228 &	2190791538023 \\
3739 &	11 &	0.00021638831563001181 &	281259932931 \\
3779 &	31 &	8.11514180382761610358 &	29319627255957776379402310856 \\
3803 &	15 &	617.93079127720341400907 &	6125701430298575603 \\
3851 &	25 &	0.60046598691752433358 &	1242978436326110671361250 \\
3907 &	7 &	0.03397552970283538243 &	1708500805 \\
3923 &	23 &	0.06601853781292555224 &	16160511049152622462466 \\
3931 &	11 &	0.00034265756041784101 &	281879748512 \\
3947 &	17 &	1.75074743805447665390 &	28850358440768828838 \\
4003 &	13 &	0.00530457688125282635 &	13618654340552 \\
4019 &	19 &	14.89793441535627886639 &	10201728430786087533096 \\
4027 &	9 &	0.06278002336060817839 &	27507874104 \\
4051 &	11 &	0.00000004836198044776 &	32711170970 \\
4091 &	33 &	0.00015048845055910840 &	13119103199363816872272275835 \\
4099 &	15 &	0.00009653216942583305 &	135256029258383 \\
4139 &	19 &	0.00018105359839238306 &	80074451508066303468 \\
4211 &	23 &	0.00001204796805032161 &	4490981695211565226634 \\
4219 &	15 &	0.03461154251357343869 &	12962166563665575 \\
4243 &	9 &	0.00066320718752067075 &	4758420837 \\
4259 &	35 &	0.32159852619627858305 &	16759504772930116391482641893012 \\
4283 &	21 &	8773.29312623157893512197 &	1535643789741168641279699 \\
4339 &	17 &	0.00136435844805308578 &	41415425722868069 \\
4363 &	9 &	0.00169672983900952063 &	12072493829 \\
4451 &	29 &	0.52221992024973798611 &	57243468904104560776219171927 \\
4483 &	9 &	345.39066341572129604016 &	90594797052049 \\
4507 &	13 &	50.23492940730687592262 &	8000361069297882 \\
4523 &	21 &	0.47925191687294653824 &	12696689618269479222463 \\
4547 &	17 &	0.00237327326500158548 &	3365011327793496632 \\
4603 &	7 &	0.67370851258063680284 &	37507414772 \\
4643 &	13 &	0.05750243316638869235 &	143101155917879067 \\
4651 &	17 &	0.00001544101918254929 &	15235214493740502 \\
4691 &	21 &	0.00000333535197828509 &	2621016573927853606584 \\
4723 &	9 &	0.95386187065080598582 &	5224661311515 \\
4787 &	25 &	196299950.09431612703528249287 &	1407786352476914766515861602276 \\
4931 &	35 &	0.00000002277116284325 &	448007986303112905715567652356 \\
4987 &	9 &	0.08774530421699530811 &	1780109232599 \\
5003 &	15 &	10706.38070422879767696786 &	34395116945398635220139 \\
5011 &	21 &	0.00050859723642713633 &	91159910603652230720 \\
5051 &	29 &	8.49904519468908134660 &	1141170772780262770088006751376 \\
5059 &	19 &	0.00156551345932179882 &	9867629103745428509 \\
5099 &	39 &	0.20398131340164447660 &	1184616873402209140039591463547587231 \\
5107 &	7 &	0.41759384988925603185 &	40655661145 \\
5147 &	19 &	9.89685243132721839321 &	110480338090747788980762 \\
5171 &	35 &	0.36860175091251048590 &	7458685927940450075817149649785164 \\
5179 &	11 &	0.00241753460843546585 &	432857713333757 \\
5227 &	15 &	0.02590550552009961454 &	28726219588915216 \\
5323 &	15 &	0.63466138387854294475 &	101661479910418479 \\
5347 &	13 &	48.38774889967386398841 &	60377540719495011 \\
5387 &	23 &	46969768.16794671957479225493 &	59861891852089030100409766615 \\
5419 &	13 &	0.00155877501356795214 &	3346444025253455 \\
5443 &	9 &	0.00772315292672455124 &	259851370654 \\
5483 &	17 &	0.00000850543167234778 &	12838230594323889723 \\
5507 &	23 &	1196283.61113816912742549060 &	308977010562437453037385949809 \\
5531 &	23 &	0.00787997073850461781 &	114928716348176366112944426 \\
5563 &	15 &	0.17610681503569755604 &	289650044515887808 \\
5651 &	31 &	0.00001430611962530720 &	1053842731936249031364945295419 \\
5659 &	19 &	0.00000095193065835809 &	2612065297147014949 \\
5683 &	11 &	44.34550258638220009604 &	12037895727172693 \\
5779 &	13 &	0.00020942032710240135 &	1281689314164362 \\
5827 &	15 &	0.00301443397675326391 &	76300991821550264 \\
5843 &	25 &	5443491913.33884952923429281132 &	965194290530223430718350605440606 \\
5851 &	21 &	0.00045098763594953089 &	4466493231539670036837 \\
5867 &	21 &	2.55857313302494265980 &	6273825474355118152363875 \\
5923 &	7 &	0.58161114692979764321 &	287690477472 \\
5939 &	35 &	603122.99403995879567858207 &	1661965438061076810127871901025304745853 \\
5987 &	15 &	110.46116247014125681999 &	13865640754944135873431 \\
6011 &	27 &	4.86492240549697162717 &	87216830057513280930236553759140 \\
6043 &	9 &	0.00721258251692584701 &	635377400757 \\
6067 &	15 &	0.00083214765715847858 &	66156492166979308 \\
6091 &	15 &	0.00004657801554426414 &	37563167157316048 \\
6131 &	31 &	0.47402035862479912750 &	10588808624066473635049652220685900 \\
6163 &	11 &	0.00109962773126547504 &	212505143679808 \\
6203 &	17 &	16956.91608309776740115902 &	17758388633281304578194836 \\
6211 &	15 &	0.00006166346655870246 &	143325584038005751 \\
6299 &	43 &	2.05775857985552601858 &	1388856212566707742241274725107204020014609 \\
6323 &	21 &	0.15649673018585925176 &	5097012104044281804792939 \\
6379 &	17 &	0.00000844841271024925 &	2826990910532861162 \\
6427 &	9 &	0.04949970413006985692 &	3642618058167 \\
6451 &	17 &	0.00000000015961944529 &	4513598851205200 \\
6491 &	31 &	0.00000661570599992748 &	26061044582689107840212347433528 \\
6547 &	11 &	0.00049212375667197424 &	28027440443663 \\
6563 &	23 &	331287.05252863171559172025 &	875165582628365469003240470630 \\
6571 &	15 &	0.00113860246370370679 &	6769437817419744845 \\
6619 &	13 &	0.01464489531673771186 &	228222213173218843 \\
6659 &	23 &	0.00004608564886020156 &	3960701905118801624128927055 \\
6691 &	21 &	0.00076792513778405932 &	44437713772698490604378 \\
6763 &	9 &	0.02842453885205402194 &	12452915791145 \\
6779 &	39 &	4.86957703804081099214 &	69700200603218027480784056371040211086176 \\
6803 &	19 &	3887.42979178406639019887 &	18785753293448529197999195649 \\
6827 &	17 &	14.01183337740201353984 &	3988442850286746357347994 \\
6883 &	9 &	0.00267729911739092469 &	8323175699073 \\
6899 &	35 &	0.00001267331994634709 &	216639728383074935163216743914068125 \\
6907 &	17 &	0.00034176219626136351 &	7191713364116285807 \\
6947 &	29 &	9807687.61625312683632853546 &	1396721892374430002050625825823542808 \\
6971 &	45 &	144.95259877715406861536 &	16300433495176531353395464824906912329257404117 \\
7019 &	43 &	0.00006851095988046843 &	419823130126477043502027992927631092428380 \\
7027 &	11 &	0.00126334189165518431 &	372846808979421 \\
7043 &	23 &	79.57658531732701207353 &	42967244325677025460819829585 \\
7187 &	25 &	1941.80715024988578872060 &	111790555942583682538467190725164 \\
7211 &	35 &	0.00000013542177935543 &	92297726998496877729409205982930880 \\
7219 &	15 &	0.00005574430738945572 &	1293006824668014040 \\
7243 &	13 &	12.33651259315172930941 &	1756240279022402156 \\
7283 &	25 &	2.49002172917008268374 &	10685692353679063484178785461948 \\
7307 &	25 &	62171.26156087145465015128 &	1419081250825229997211797415503825 \\
7331 &	33 &	0.00000000032502361559 &	219455139795728684736426909205613 \\
7411 &	25 &	0.00015095236624874787 &	283024230782339795658408815 \\
7451 &	35 &	93950.78490910382748157143 &	1553319522388195483106386505386846710579688 \\
7459 &	15 &	0.00003036649237987241 &	2337968558206904465 \\
7499 &	33 &	0.00000135124978504196 &	60792334311439101075808649223338741 \\
7507 &	11 &	0.00010611225020431393 &	506903678439547 \\
7523 &	35 &	33522.02988333545407545411 &	15470691108260141377130066889443447728251 \\
7547 &	15 &	0.60032452144273867428 &	258920923898858718493581 \\
7603 &	11 &	0.00653892833763582491 &	1434870589808453 \\
7643 &	29 &	179270.99125915725948310242 &	2288294534095276709215678772272550333 \\
7691 &	43 &	1936270.78489178696965784529 &	3624801748457692072742156116148403862768295124739 \\
7699 &	27 &	0.00000699702201774193 &	2304258410337272194320291613 \\
7723 &	9 &	0.01895970678708005242 &	46425972727345 \\
7867 &	11 &	1.61910857940903382775 &	969216189901440388 \\
7883 &	17 &	0.23552172495493403317 &	28080473969786372085402707 \\
7907 &	21 &	0.00175081994226267973 &	568887022906905104307160000 \\
7963 &	13 &	0.00002874140087673655 &	1773163400695878 \\
8011 &	25 &	0.00000000000000681929 &	5972318120577118158596 \\
8059 &	21 &	0.00000000008523807585 &	384252254712642261955 \\
8123 &	21 &	2383.30719868549773975742 &	18940429636369330011648211950850 \\
8147 &	37 &	12176247.61615672604210432982 &	108543217201536511344104411121816314211436357 \\
8171 &	21 &	0.00001153050501312074 &	43709209245428902074487555857 \\
8179 &	25 &	0.00000000000408693701 &	245379074777785127542855 \\
8219 &	35 &	0.00003219371470307854 &	320575044136527674847638735156405728284 \\
8243 &	21 &	0.00003184119831052584 &	169936812859198942971657899 \\
8291 &	47 &	13.82276752176918097863 &	144216236301930044357917532496691060210143449492929 \\
8363 &	35 &	232.06334304127124142727 &	15312049368913584760165020656408258883985 \\
8387 &	21 &	1545.33342862710955449647 &	646162618485014801463564772032 \\
8419 &	19 &	0.00000012646740535682 &	7680287320591089475912 \\
8443 &	11 &	1.17893537822797473748 &	35478826850611232 \\
8467 &	15 &	0.00021223474007323571 &	1265829019458818985 \\
8539 &	17 &	0.00000478356805956157 &	2275177608844149620235 \\
8563 &	9 &	2.03683466086338426663 &	5420212257148597 \\
8627 &	21 &	3853.12162724028605983782 &	198171648305412280420765428378075 \\
8699 &	35 &	0.00002481823232900755 &	11240439275853225767486298164708391703500 \\
8707 &	15 &	0.05248611135357774651 &	544957212659766295149 \\
8731 &	17 &	0.00003500145787893966 &	1576791412901313685993 \\
8747 &	21 &	484018.07264403016140460258 &	288450573402627146790548293522904 \\
8779 &	15 &	0.00000030582300865004 &	3039062315648643705 \\
8803 &	9 &	0.01552030193111727803 &	126149700575776 \\
8819 &	49 &	4034780.63025405890013810127 &	23921217585664253728980536953688630289914216264814750750 \\
8867 &	27 &	2.17002992458273671577 &	128281233388318287503394675078625248 \\
8923 &	19 &	159.33124314583988814840 &	14174087757188097620330473 \\
8963 &	29 &	1985418.83735640580162530103 &	1118129255137682166151656920634902507171 \\
8971 &	19 &	0.00000000983100532208 &	2258409802656205931201 \\
9011 &	33 &	0.00335393090412302635 &	3045909543390800207752459587622105801155 \\
9043 &	15 &	0.10664966467865754486 &	247161131467830177656 \\
9059 &	39 &	13.77046293947696022352 &	1591196881664881653149442733268746342271235079 \\
9067 &	9 &	0.00731970822333483793 &	111048528074239 \\
9091 &	21 &	0.00000000236926748020 &	236934738121647151313216 \\
9187 &	21 &	32.52092107710830419611 &	366682928077533301738876375 \\
9203 &	31 &	1023606763.23245034828055208642 &	10450577199144964797897931682954499342069972 \\
9227 &	25 &	0.03059726342701622948 &	61065081813711319346093946656236 \\
9283 &	11 &	0.00312668644337035262 &	114634424721388682 \\
9323 &	29 &	158279113.16389204676979182791 &	198428451440095157640040934380936648153653 \\
9371 &	49 &	0.00000000381610952189 &	9868042349887290226589618012608800982902276255554 \\
9403 &	11 &	0.00062444568609082732 &	9557020329789716 \\
9419 &	35 &	0.00004080832458224277 &	7755628001204338180900900374435880983024 \\
9467 &	41 &	85262.97994248194021607026 &	712735551990570229824239164893859490851236314718 \\
9491 &	45 &	1.93118062921922617356 &	699904157096559596488802074953836695863315457589115 \\
9539 &	55 &	0.00000011262029208767 &	1854017230259878781685683081695866642830423627047547360 \\
9547 &	13 &	4.55503594583652426666 &	636512933633007687775 \\
9587 &	23 &	100123.63380876703245300853 &	36563316420884623249966367262977278 \\
9619 &	19 &	0.00000012220071093949 &	5176494409557611828238 \\
9643 &	11 &	0.00824565865962835292 &	516445608289196035 \\
9739 &	13 &	0.00000131331766917028 &	3753135247337182744 \\
9787 &	11 &	0.26005393619869399226 &	362054016521593344 \\
9803 &	37 &	5930542264.64370531353643671288 &	3056997114849536281991895391097536748139790872620 \\
9811 &	21 &	0.00000000057298287826 &	127050303874152535296768 \\
9851 &	45 &	0.00000000000000654205 &	18104457473462858976938437755791606152672401 \\
9859 &	21 &	0.00000242965527955730 &	28849640043711849588277760 \\
9883 &	17 &	3.13476671971555114510 &	242746734273754376709691 \\
9907 &	15 &	0.00104023793825974741 &	106327569146508318840 \\
9923 &	25 &	0.00000242632276300466 &	2387019558186614148757292853595 \\
9931 &	23 &	0.00000000000063190232 &	807612994249556860279774 \\
10067 &	21 &	518.08680600795693371589 &	342288622319158792879184584048129 \\
10091 &	39 &	0.00000478306284603539 &	33197255379943602421429579685905522266248888 \\
10099 &	25 &	0.00000000001780723538 &	52822827121201782525001225 \\
10139 &	55 &	73803.78400959271904263740 &	33150521228939817885342694776696740143379123299513891375082064 \\
10163 &	39 &	67991.83407099603855149887 &	440513245214192270403084087547571959843972111425 \\
10211 &	43 &	36.21163226407725524787 &	2018436122269325072524514224058369879912299468812420 \\
10243 &	15 &	0.08816299542047442521 &	1839984519589335519209 \\
10259 &	43 &	652.29068417959804070259 &	718692371932632451954249149920382743582780046485373 \\
10267 &	15 &	0.00014937068406596087 &	382939450567034881777 \\
10331 &	55 &	0.00062202944015487748 &	6664989351016927856864678973892159602391574133103405618788 \\
10427 &	31 &	43.44997691536163284782 &	2845156404392349071355292880344541653644 \\
10459 &	15 &	0.00028151807209305816 &	4605303213234756184743 \\
10499 &	41 &	507989.87585259730539676742 &	36955062681702947421893940479768814784946674217168545 \\
10531 &	27 &	0.00000000000146797238 &	3351769586174993494323439275 \\
10627 &	9 &	0.00017318804748393473 &	398694874516979 \\
10651 &	15 &	0.00011268583573677928 &	12553194111400362264301 \\
10667 &	39 &	25506036465736.67223882308179488976 &	44811079267440631723932527212656541784290245322216891 \\
10691 &	45 &	3520.56248376976938840747 &	232694011311673350014283799021469164888464834922982077 \\
10723 &	15 &	0.00000092340639534283 &	21063173152441592413 \\
10739 &	37 &	0.00000000004343365251 &	63654774159741723864112646151198426204132 \\
10771 &	21 &	0.00006625182262625102 &	6376645704546135373491729472 \\
10859 &	45 &	0.00001929088385542839 &	144220156391029884341565008365356176802926959516032 \\
10867 &	23 &	0.00238767061834684639 &	10525038031112725310622398635 \\
10883 &	23 &	515416.17141491515842049134 &	87346650878073993681932895672040241977 \\
10891 &	19 &	0.00000016242472588156 &	11884817420032518255654205 \\
10939 &	27 &	0.00000000901518555777 &	384273411889144154980660519291 \\
10979 &	35 &	0.00368603736420883092 &	3439622056704005260743458738409400370601464092 \\
10987 &	11 &	0.00349823665853601300 &	20844399361215150 \\
11003 &	27 &	2285808.76686708228105121178 &	215352519101955536146438795127413473169151 \\
11027 &	27 &	1.79342141479218417410 &	1218295757801391198254967804816697855 \\
11059 &	25 &	0.00000017029649146888 &	339791808371879703974813638561 \\
11083 &	15 &	0.00029395864104112152 &	69130942815185383643 \\
11131 &	25 &	0.00000000011466927704 &	779081246955355170396973552 \\

\end{longtable}
\end{center}

After all the calculations were finished, John Coates informed us that for $q$ congruent to $3$ modulo $8$,  
one should actually take the square root of the $j(\mathcal O_K) - 12^3$ with negative imaginary part to get the 
appropriate quadratic twist of the Gross curve $A$ (see the formula (2.2)  of  the paper by Buhler and Gross \cite{BG}).  
Hence, we should work with the conjugate of the equation (1) under the Galois group of $H$ over $\Bbb Q$.  
We have checked for a few small $q$ congruent to $3$ modulo $8$, that the $L$-values, torsion parts and 
Tamagawa numbers of the two conjugate curves are the same. As a consequence, the analytic orders of 
Tate-Shafarevich groups of these curves are the same. John Coates expects that Tate-Shafarevich groups 
of these curves are actually isomorphic, but all is not totally clear theoretically.

\bigskip 

{\bf Acknowledgements.} 
We primary thank John Coates for suggesting the problem to us, and for very inspiring correspondence.  
We would like to thank Mark Watkins \cite{Wat}, who sent us a new algorithm allowing to make great progress 
in our  calculations of the values $L(E/H,1)$.

Institute of Mathematics, University of Szczecin, Wielkopolska 15, 
70-451 Szczecin, Poland; E-mail addresses: andrzej.dabrowski@usz.edu.pl and dabrowskiandrzej7@gmail.com;  
tjedrzejak@gmail.com;   lucjansz@gmail.com

\end{document}